\documentclass[11pt]{article}
\usepackage{amsmath, amssymb, amsthm}
\usepackage{geometry}
\usepackage{hyperref}
\hypersetup{
    colorlinks=true,
    linkcolor=blue,
    citecolor=blue,
    urlcolor=blue
}

\newtheorem{theorem}{Theorem}[section]  
\newtheorem{lemma}[theorem]{Lemma}      

\newtheorem{proposition}[theorem]{Proposition}
\theoremstyle{definition}
\newtheorem{definition}[theorem]{Definition}

\theoremstyle{remark}

\newcommand{\keywords}[1]{%
  \vspace{0.5cm}
  \noindent\textbf{Keywords:} #1
  \vspace{0.5cm} }
\title{\textbf{On the Explicit Expression of an Extended Version of Riemann Zeta Function}}
\author{\textbf{Author}: Yushi Huang, Jiangsu Tianyi High School \\ \newline \textbf{Adviser}: Ting Wu, Department of Mathematics, Nanjing University}
\date{}
\begin{document}
\maketitle
\begin{abstract}
In this paper, we focus on the explicit expression of an extended version of Riemann zeta function. We use two different methods, Mellin inversion formula and Cauchy's residue theorem, to calculate a Mellin-Barnes type integral of the analytic function regarding $z$: $\Gamma(z)\Gamma(s-z)u^{-z}$ ($u\in (0,1)$, $s\in \mathbb{C}$). We provide the necessary background on the analytic properties of Gamma and Riemann zeta function to confirm the absolute convergence of this Mellin-Barnes integral. Next, we represent the extended version of Riemann zeta function $\sum_{m=1}^{\infty}\sum_{n=1}^{\infty}{(m+n)^{-s}}$ using the following complex integral where the real part of $s$ is larger than 2 and $c>1$ is chosen to make $\Re(s)-c$ larger than 1.  
$$\Gamma(s)\sum_{m=1}^{\infty}\sum_{n=1}^{\infty}{(m+n)^{-s}}=\frac{1}{2\pi i} \int_{c - i\infty}^{c + i\infty} \zeta(z) \zeta(s - z) \Gamma(z) \Gamma(s - z) \, dz$$
We provide the evaluation of this integral by changing the integration path from straight line $\Re(z)=c$ into a rectangular contour whose left side is positioned at negative infinity. We apply the functional equation of Riemann zeta function, Euler's reflection formula, and Legendre's duplication formula to evaluate the integral segment through $\Re(z)=-\infty$. After introducing Hurwitz zeta function and properly calculating the difference between the sum of residues in two analogous rectangular contours, we finalize the evaluation. Lastly, we demonstrate the connection of this result with other intricate integrals involving special functions, such as the hyperbolic function. Additionally, we discuss its applications in deriving explicit expressions for the Barnes zeta function. 
\end{abstract}
\newpage
\tableofcontents
\newpage
\section{Introduction}
\paragraph{}
Gamma function $\Gamma(s)$, which was originally introduced by Leonhard Euler in 1729 in his letter to Goldbach \cite{davis1959euler}, is defined as an infinite product: for all 
\( s \in \mathbb{C} \), 
\begin{equation}
\frac{1}{\Gamma(s)}=s\prod_{n=1}^{\infty}\left(1+\frac{1}{n}\right)^{-s}\left(1+\frac{s}{n}\right)
\label{eqgamma}
\end{equation}   

Later, Daniel Bernoulli derived an integral representation of Gamma function for all complex number $s$ where the real part of $s$ is larger than 0: 
$$\Gamma(s)=\int_{0}^{\infty} t^{s-1}e^{-t} dt$$

Gamma function is then analytically continued to a meromorphic function which is holomorphic over the entire complex plane except for its simple poles at zero and negative integers. It plays an important role in mathematical analysis for calculating factorials and serves as a powerful tool in solving intricate integrals. Likewise, Riemann zeta function, which is equally important in modern mathematics especially analytic number theory, is originally defined by Lenohard Euler in 1737 as: 
$$\zeta(s)=\sum_{n=1}^{\infty}\frac{1}{n^s}$$
where $s$ is a real number. Subsequently, Bernhard Riemann extended $s$ in Euler's definition to a complex variable \cite{Riemann1859}, proved its meromorphic continuation and functional equation, and intensively studied its deep connection with the distribution of prime numbers. Riemann zeta function and Gamma function are closely related through the following integral \cite{apostol1976introduction}:
\begin{equation}
\Gamma(s)\zeta(s)=\int_{0}^{\infty}\frac{x^{s-1}}{e^x-1}dx 
\label{cond}
\end{equation} 

Historically, mathematicians are interested in studying complex integrals involving Gamma and Riemann zeta function. In 1908 and then 1910, Ernest William Barnes studied two types of integrals involving solely a product of Gamma functions, known as the first and second Barnes lemma \cite{Barnes1908}: 
\[
\frac{1}{2\pi i} \int_{-i\infty}^{i\infty} \Gamma(a+s)\Gamma(b+s)\Gamma(c-s)\Gamma(d-s)\, ds 
\]
\[
\frac{1}{2\pi i} \int_{-i\infty}^{i\infty} \frac{\Gamma(a+s)\Gamma(b+s)\Gamma(c+s)\Gamma(1-d-s)\Gamma(-s)}{\Gamma(a+b+c-d+1+s)}\, ds
\]

In 1997, Masanori Katsurada studied complex integrals involving both Gamma and Riemann zeta function \cite{Katsurada1997}: 
\[ 
\frac{1}{2\pi i} \int_{\sigma-i\infty}^{\sigma+i\infty} \frac{\Gamma(a+s)\Gamma(-s)}{\Gamma(b+s)} \zeta(c+s) (-z)^s ds
\]

Their work has provided important explicit expressions for complicated integrals or infinite sums, establishing connections with other special functions, such as the hyperbolic function and hypergeometric function. 

In section \ref{section22} of this paper, we first derive an identity containing a Mellin-Barnes type integral: 
\begin{equation}
\Gamma(s)(1 + u)^{-s}=\frac{1}{2\pi i} \int_{c- i\infty}^{c + i\infty} \Gamma(z) \Gamma(s - z) u^{-z} \, dz 
\label{eq1234}
\end{equation}

Then, in section \ref{section23} we get an important relationship between a double sum version of Riemann zeta function and a complex integral involving both Gamma and Riemann zeta function: 
\begin{equation}
\Gamma(s)\sum_{m=1}^{\infty}\sum_{n=1}^{\infty}{(m+n)^{-s}}=\frac{1}{2\pi i} \int_{c - i\infty}^{c + i\infty} \zeta(z) \zeta(s - z) \Gamma(z) \Gamma(s - z) \, dz
\label{eq123}
\end{equation}

In section \ref{section24}, we apply the analytic properties of both Gamma and Riemann zeta function and Cauchy's residue theorem to evaluate this complex integral \eqref{eq123}. In the last section \ref{section3}, we discuss the applications of this result and its efficacy for further research. 

\section{Main manipulations} \label{section2}
\paragraph{}
In this section, we tackle the integral \eqref{eq123}.  
\subsection{Preliminaries on absolute convergence} \label{section21}
\paragraph{}
We first introduce a few analytic properties of Gamma function, including its asymptotic approximation, analytic behavior in vertical strips, and upper bound, which are necessary tools for our work in the subsequent sections. Later, we use these results to derive several propositions regarding the Mellin-Barnes type integral \eqref{eq1234}. 
\begin{lemma}
Let $0<\delta<\pi$. Then for $-\pi+\delta<\arg s<\pi-\delta$,
\begin{equation}
\log \Gamma(s)=\left(s-\frac{1}{2}\right) \log s-s+\log \sqrt{2 \pi}+\mathcal{O}\left(|s|^{-1}\right), \quad|s| \rightarrow \infty \label{con:1}
\end{equation}
\end{lemma}

\begin{proof}
By Euler's definition of Gamma function \eqref{eqgamma}, for $-\pi<\arg s<\pi$ we have: 
\begin{equation}
\begin{aligned}
\log \Gamma(s) & =-\log s+\sum_{n=1}^{\infty}\left\{s \log \left(1+\frac{1}{n}\right)-\log \left(1+\frac{s}{n}\right)\right\} \\
& =-\log s+\lim _{N \rightarrow \infty}\left\{s \log N-\sum_{n=1}^{N} \log \left(1+\frac{s}{n}\right)\right\} 
\end{aligned}\label{con:2}
\end{equation}

By partial summation,
\begin{equation}
\sum_{n=1}^{N} \log \left(1+\frac{s}{n}\right)=\log (1+s)+\int_{1}^{N} \log \left(1+\frac{s}{u}\right) d[u] \label{con:3}
\end{equation}

Note that $[u]=u-\frac{1}{2}+\rho(u)$ where $\rho(u)=1 / 2-\{u\}$. Thus the above integral becomes
\[
\int_{1}^{N} \log \left(1+\frac{s}{u}\right) d u+\int_{1}^{N} \log \left(1+\frac{s}{u}\right) d \rho(u)
\]

We have: 
\[
\begin{aligned}
\int_{1}^{N} \log \left(1+\frac{s}{u}\right) d u & =N \log \left(1+\frac{s}{N}\right)-\log (1+s)+\int_{1}^{N} \frac{s}{u+s} d u \\
& =N \log \left(1+\frac{s}{N}\right)-(s+1) \log (1+s)+s \log (N+s)
\end{aligned}
\]
and
\[
\int_{1}^{N} \log \left(1+\frac{s}{u}\right) d \rho(u)=\frac{1}{2} \log \left(1+\frac{s}{N}\right)-\frac{1}{2} \log (1+s)-\int_{1}^{N} \rho(u)\left(\frac{1}{u+s}-\frac{1}{u}\right) d u
\]

Thus, we get: 
\[
\begin{aligned}
\sum_{n=1}^{N} \log \left(1+\frac{s}{n}\right)= & -\left(s+\frac{1}{2}\right) \log (1+s)+\left(N+s+\frac{1}{2}\right) \log \left(1+\frac{s}{N}\right) \\
& +s \log N-\int_{1}^{N} \rho(u)\left(\frac{1}{u+s}-\frac{1}{u}\right) d u
\end{aligned}
\]

Back to (\ref{con:2}), we get: 
\[
\begin{aligned}
\log \Gamma(s)= & \left(s+\frac{1}{2}\right) \log (1+s)-\log s+\int_{1}^{\infty} \frac{\rho(u)}{u+s} d u-\int_{1}^{\infty} \frac{\rho(u)}{u} d u \\
& -\lim _{N \rightarrow \infty}\left(N+s+\frac{1}{2}\right) \log \left(1+\frac{s}{N}\right)
\end{aligned}
\]

The last limit is equal to $s$. Thus,
\begin{equation}
\log \Gamma(s)=\left(s+\frac{1}{2}\right) \log (1+s)-\log s-s+\int_{1}^{\infty} \frac{\rho(u)}{u+s} d u+c_{0}\label{con:4}
\end{equation}
we write
\[
\sigma(x)=\int_{1}^{x} \rho(u) d u=\frac{\{x\}}{2}-\frac{\{x\}^{2}}{2}
\]

Then
\begin{equation}
\int_{1}^{\infty} \frac{\rho(u)}{u+s} d u=\int_{1}^{\infty} \frac{\sigma(u)}{(u+s)^{2}} d u \ll \int_{1}^{\infty} \frac{d u}{|u+s|^{2}} \label{con:5}
\end{equation}

One has
\[
\left|\int_{1}^{\infty} \frac{d u}{|u+s|^{2}}\right|\leq \int_{1}^{\infty} \frac{d u}{u^{2}+|s|^{2}-2 u|s||\cos \delta|}=\frac{1}{|s|} \int_{\frac{1}{|s|}}^{\infty} \frac{d t}{1+t^{2}-2 t|\cos \delta|}.
\]

Note that: 
\[
1+t^{2}-2 t|\cos \delta|=(1-t|\cos \delta|)^{2}+\left(1-|\cos \delta|^{2}\right) t^{2}
\]

Hence for $|s|>1$,
\begin{equation}
\int_{1}^{\infty} \frac{d u}{|u+s|^{2}} \ll \frac{1}{|s|} \int_{\frac{1}{|s|}}^{1} \frac{d t}{(1-t|\cos \delta|)^{2}}+\frac{1}{|s|} \int_{1}^{\infty} \frac{d t}{\left(1-|\cos \delta|^{2}\right) t^{2}} \ll_{\delta} \frac{1}{|s|} \label{con:6}
\end{equation}

This shows that
\[
\int_{1}^{\infty} \frac{\rho(u)}{u+s} d u=\mathcal{O}_{\delta}\left(\frac{1}{|s|}\right)
\]

Since
\[
\log (1+s)=\log s+\frac{1}{s}+\mathcal{O}_{\delta}\left(\frac{1}{|s|^{2}}\right), \quad|s| \rightarrow \infty
\]
we get from (\ref{con:4})--(\ref{con:6}) that
\[
\log \Gamma(s)=\left(s-\frac{1}{2}\right) \log s-s+c_{0}+1+\mathcal{O}_{\delta}\left(|s|^{-1}\right)
\]

Taking $s=n, n+1 / 2,2 n$ and using the Legendre's duplication formula we get $c_{0}=\log\sqrt{2 \pi}$.
\end{proof}

\begin{lemma}
For $\alpha \leq \sigma \leq \beta$, there holds
\[
\Gamma(\sigma+i t)=\sqrt{2 \pi}|t|^{\sigma+i t-\frac{1}{2}} e^{-\frac{\pi|t|}{2}-i t+i \frac{\pi}{2}\left(\sigma-\frac{1}{2}\right) \operatorname{sgn}(t)}\left(1+\mathcal{O}\left(|t|^{-1}\right)\right), \quad|t| \rightarrow+\infty
\]
where the implied constant depends on $\alpha$ and $\beta$.
\label{lemma22}
\end{lemma}

\begin{proof}
By (\ref{con:1}) we get    
\[
\Gamma(s)=\sqrt{2 \pi} e^{\left(s-\frac{1}{2}\right) \log s-s}\left(1+\mathcal{O}\left(|s|^{-1}\right)\right)
\]

Let $s=\sigma+i t$. For $|t|>|\sigma|$ we have
\[
\log s=\log i t+\log \left(1+\frac{\sigma}{i t}\right)=\frac{\pi i \operatorname{sgn}(t)}{2}+\log |t|+\frac{\sigma}{i t}+\mathcal{O}\left(|t|^{-2}\right)
\]
where the $\mathcal{O}$-constant depends on $\alpha$ and $\beta$. Thus
\[
\begin{aligned}
\left(s-\frac{1}{2}\right) \log s-s& =\frac{\pi i \operatorname{sgn}(t)}{2}\left(\sigma+i t-\frac{1}{2}\right)+\left(\sigma+i t-\frac{1}{2}\right) \log |t|+\left(\sigma+i t-\frac{1}{2}\right) \frac{\sigma}{i t}-\sigma-i t+\mathcal{O}\left(|t|^{-1}\right) \\
& =\frac{\pi i \operatorname{sgn}(t)}{2}\left(\sigma-\frac{1}{2}\right)-\frac{\pi|t|}{2}+\left(\sigma+i t-\frac{1}{2}\right) \log |t|-i t+\mathcal{O}\left(|t|^{-1}\right)
\end{aligned}
\]

Therefore, 
\begin{equation}
\Gamma(\sigma+i t)=\sqrt{2 \pi}|t|^{\sigma+i t-\frac{1}{2}} e^{-\frac{\pi|t|}{2}-i t+i \frac{\pi}{2}\left(\sigma-\frac{1}{2}\right) \operatorname{sgn}(t)}\left(1+\mathcal{O}\left(|t|^{-1}\right)\right)\label{con:7}
\end{equation}
\end{proof}
\begin{lemma}
    For all \( \sigma \geq \frac{1}{2} \), there is a constant \( C(\sigma) \) such that
\[
|\Gamma(x + iy)| \leq C(\sigma) e^{-|y|} \ \ \ \ \left(\frac{1}{2} \leq x \leq \sigma, \, y \in \mathbb{R}\right)
\]\label{con:8}
\end{lemma}
\begin{proof}
By (\ref{con:7}) we get:
\[
\begin{aligned}
\Gamma(x + iy)=&\sqrt{2 \pi}|y|^{x+i y-\frac{1}{2}} e^{-\frac{\pi|y|}{2}-i y+i \frac{\pi}{2}\left(x-\frac{1}{2}\right) \operatorname{sgn}(y)}\left(1+\mathcal{O}\left(|y|^{-1}\right)\right)\\
=&\sqrt{2 \pi}|y|^{x-\frac{1}{2}}e^{-\frac{\pi|y|}{2}}e^{i \frac{\pi}{2}\left(x-\frac{1}{2}\right) \operatorname{sgn}(y)}\left(\frac{|y|}{e}\right)^{iy}\left(1+\mathcal{O}\left(|y|^{-1}\right)\right)
\end{aligned}
\]

Hence, 
\[
|\Gamma(x + iy)|=\sqrt{2 \pi}e^{-\frac{\pi|y|}{2}}|y|^{x-\frac{1}{2}}|\left(1+C_0(\sigma)|y|^{-1}\right)|
\]

If $|y|\geq 1$ and $x\leq\frac{3}{2}$, then we have:  
\[
|\Gamma(x + iy)|\leq C(\sigma) e^{-|y|}
\]
where the $C(\sigma)=\sqrt{2 \pi}\left(1+|C_0(\sigma)|\right)$.

If $|y|\geq 1$ and $x\geq\frac{3}{2}$, then we have: 
\[
|y|^{x-\frac{1}{2}}\leq|y|^{\sigma-\frac{1}{2}}\leq C_1(\sigma)e^{(\frac{\pi}{2}-1)|y|}
\]

Thus, we have: 
\[
|\Gamma(x + iy)|\leq C(\sigma) e^{-|y|}
\]
where the $C(\sigma)=\sqrt{2 \pi}|\left(1+C_0(\sigma)\right)|C_1(\sigma)$.\\

If $|y|\leq 1$, then we have: 
\[
|\Gamma(x + iy)|\leq C(\sigma) e^{-|y|}
\]
where the $C(\sigma)=\sqrt{2 \pi}\left(1+|C_0(\sigma)|\right)$.
\end{proof}

\begin{lemma}
    For a positive integer $k$, we have
\[
\Gamma(z - k) \Gamma(s - z + k) = (-1)^k \Gamma(z) \Gamma(s - z) \prod_{1 \leq j \leq k} \left( 1 + \frac{s - 1}{j - z} \right).
\]\label{con:9}
\end{lemma}
\begin{proof}
We apply the functional equation of Gamma function: for $s\in \mathbb{C}$, $s\Gamma(s)=\Gamma(s+1)$. 

Hence
\[
(z-k)(z-k+1)...(z-1)\Gamma(z-k)=\Gamma(z)
\]
and 
\[
(s-z)(s-z+1)...(s-z+k-1)\Gamma(s-z)=\Gamma(s-z+k)\]

Thus
\[
\begin{aligned}
\Gamma(z - k) \Gamma(s - z + k) &= \frac{\Gamma(z)}{(z-k)(z-k+1)\cdots(z-1)}\cdot (s-z)(s-z+1)\cdots(s-z+k-1)\Gamma(s-z) \\
&= \Gamma(z)\Gamma(s-z)\cdot (-1)^k \frac{(1-z+s-1)(2-z+s-1)\cdots(k-z+s-1)}{(k-z)((k-1)-z)\cdots(1-z)} \\
&= \Gamma(z)\Gamma(s-z)\cdot (-1)^k \prod_{1\le j\le k} \frac{j-z+s-1}{j-z} \\
&= (-1)^k \Gamma(z) \Gamma(s - z) \prod_{1 \leq j \leq k} \left( 1 + \frac{s - 1}{j - z} \right)
\end{aligned}
\]
\end{proof}
\begin{lemma}
Given $s=\sigma+it$ satisfying $ \sigma = \Re (s) > 2$, if $\Re(z)=c\ge \frac{1}{2}$ and $\Re(s)-c\ge \frac{1}{2}$, there is a constant $ M(s) > 0 $ such that
\[
|\Gamma(z-k)\Gamma(s-z+k)| \leq M(s) k^{\sigma-1} e^{-|y|}\ \ \ \left(k \geq 1, y=\Im(z) \right)
\]
\label{lemmagood}
\end{lemma}

\begin{proof} 
By lemma \ref{con:8} and lemma \ref{con:9} we get:
\begin{align}
|\Gamma(z-k)\Gamma(s-z+k)| &= |\Gamma(z)||\Gamma(s-z)|\prod_{1 \leq j \leq k} \left|1 + \frac{s - 1}{j - z}\right| \nonumber\\
&\leq C_1 e^{-|y|}C_2(s)e^{-|t-y|}\prod_{1 \leq j \leq k} \left|1 + \frac{s - 1}{j - z}\right|  \label{eq8}
\end{align}

Due to the identity $\Gamma(\overline{z})=\overline{\Gamma(z)}$ and $|z|=|\overline{z}|$, the signs of $\Im(z)$ or $\Im(s)$ do not impact the value of the modulus of Gamma function. Hence, we can assume both $t$ and $y$ are larger than 0. 

We know that
\[
|1+w|\leq e^{\Re(w)+\mathcal{O}(|w|^2)}
\]
so we have: 
\[
\left|1+\frac{s - 1}{j - z}\right| \leq e^{\Re(\frac{s-1}{j-z})+\mathcal{O}(|\frac{s-1}{j-z}|^2)}
\]

We calculate the real part of $\frac{s-1}{j-z}$ and estimate the error term: 
\[
\begin{aligned}
\Re\left(\frac{s-1}{j-z}\right)&=\frac{(\sigma-1)+it}{(j-c)-iy}\\
&=\frac{(\sigma-1)(j-c)-ty}{(j-c)^2+y^2}\\
&\le \frac{(\sigma-1)(j-c)}{(j-c)^2+y^2}\\
&\le \frac{\sigma-1}{j-c}
\end{aligned}
\]
and
\begin{align*}
\mathcal{O}\left(\left|\frac{s-1}{j-z}\right|^2\right)&=\mathcal{O}\left(\frac{\left|(\sigma-1)+it\right|^2}{\left|(j-c)-iy\right|^2}\right)\\
&=\mathcal{O}\left(\frac{(\sigma-1)^2+t^2}{(j-c)^2}\right)\\
&=M_{0}(s)
\end{align*}
which is a constant only dependent on $s$. 

Hence, we can get: 
\begin{align*}
\prod_{j=1}^k \left|1+\frac{s-1}{j-z}\right| &\le \exp \left\{\sum_{j=1}^k \Re\left(\frac{s-1}{j-z}\right)+\mathcal{O}\left( \sum_{j=1}^{k} \left|\frac{s-1}{j-z}\right|^2 \right) \right\}\\
&\le \exp\left\{   \sum_{j=1}^k \frac{\sigma-1}{j-c} +M_0(s) \right\} \\
&\le \exp\{(\sigma-1)\log(k)+M_1(s)\}
\end{align*}
for some constant $M_1(s)$ dependent on $s$. Here, we use the upper bound for harmonic number $H_k=\sum_{j=1}^k \frac{1}{j}\le \log(k)+1$. 
Thus, from \eqref{eq8}, we get: 
\begin{align*}
|\Gamma(z-k)\Gamma(s-z+k)|&\le C_1 e^{-|y|} C_2(s) e^{-|t-y|} e^{(\sigma-1)\log(k)+M_1(s)} \\
&\le C_1 C_2(s) e^{M_1(s)} k^{\sigma-1} e^{-|y|}\\
&=M(s) k^{\sigma-1} e^{-|y|}
\end{align*}
where $M(s)=C_1 C_2(s) e^{M_1(s)}$. 
\end{proof}

So far, we have finished the properties related to Gamma function, we will then deduce two propositions about the Mellin-Barnes type integral \eqref{eq1234} we are concerned about. The first one validates the absolute convergence of this type of integral, and the second evaluates its value after shifting the integration path to negative infinity. 
\begin{proposition}
    For \( s \in \mathbb{C}, \ \Re(s)>2, \ u \in (0,1), \ k \geq 1 \), $c\ge \frac{1}{2}$, and $\Re(s)-c\ge \frac{1}{2}$, the following complex integral
\[
\frac{1}{2\pi i} \int_{c - k - i \infty}^{c - k + i \infty} \Gamma(z)\Gamma(s-z)u^{-z} dz
\]
converges absolutely.
\label{propsotioncool}
\end{proposition}

\begin{proof}
We substitute $z$ in the integrand by $z-k$: 
$$\frac{1}{2\pi i} \int_{c - k - i \infty}^{c - k + i \infty} \Gamma(z)\Gamma(s-z)u^{-z} dz=\frac{1}{2\pi i} \int_{c - i \infty}^{c + i \infty} \Gamma(z-k)\Gamma(s-z+k)u^{-z+k} dz$$

Use the upper bound given by lemma \ref{lemmagood}, we get: 
\begin{align*}
\left| \frac{1}{2\pi i} \int_{c - k - i \infty}^{c - k + i \infty} \Gamma(z)\Gamma(s-z)u^{-z} dz   \right| &=\left| \frac{1}{2\pi i} \int_{c - i \infty}^{c + i \infty} \Gamma(z-k)\Gamma(s-z+k)u^{-z+k} dz \right| \\
&\le \frac{1}{2\pi i} \int_{c - i \infty}^{c + i \infty} \left| \Gamma(z-k)\Gamma(s-z+k)u^{-z+k} dz \right|\\
&\le \frac{1}{2\pi } \int_{- \infty}^{+  \infty} M(s)k^{\sigma -1}e^{-|t|}u^{-c+k}dt\\
& = \frac{1}{\pi}\int_{0}^{\infty}{M(s)k^{\sigma -1}u^{-c+k} e^{-t}dt}\\
& = \frac{1}{\pi}M(s)k^{\sigma -1}u^{-c+k}
\end{align*}
which is $<\infty$. 
\end{proof}

\begin{proposition}
    Under the same condition as proposition \ref{propsotioncool}, 
     $$\lim_{k\to \infty}\frac{1}{2\pi i} \int_{c - k - i \infty}^{c - k + i \infty} \Gamma(z)\Gamma(s-z)u^{-z} dz=0$$
\label{propositioncoo}
\end{proposition}

\begin{proof}
Proposition \ref{propsotioncool} states: 
$$\left| \frac{1}{2\pi i} \int_{c - k - i \infty}^{c - k + i \infty} \Gamma(z)\Gamma(s-z)u^{-z} dz   \right|\le  \frac{1}{\pi}M(s)u^{-c} k^{\sigma -1}u^{k} $$
Since $u\in(0,1)$, $$\lim_{k\to \infty} k^{\sigma -1}u^{k} =0$$ for $\sigma=\Re(s)\ge 1$. Thus, as $k$ approaches infinity, 
$$\lim_{k\to \infty}\frac{1}{2\pi i} \int_{c - k - i \infty}^{c - k + i \infty} \Gamma(z)\Gamma(s-z)u^{-z} dz$$ tends to zero. 
\end{proof}

\subsection{Transition to complex integral} \label{section22}
\paragraph{}
Now, we introduce an important lemma that links a power function $(1+u)^{-s}$ and the Mellin-Barnes integral we have discussed in the previous section. 
\begin{lemma}
For $s\in \mathbb{C}$, $\Re(s)\ge 1$, $u,c\in (0,1)$, we have
\begin{equation}
 \Gamma(s)(1 + u)^{-s}=\frac{1}{2\pi i} \int_{c- i\infty}^{c + i\infty} \Gamma(z) \Gamma(s - z) u^{-z} \, dz\label{eq2} 
\end{equation}
\label{lemma210}
\end{lemma}
Throughout section \ref{section22}, we use two different approaches to prove the above lemma. Some of the techniques involved are paramount for our subsequent work. 

\subsubsection{First proof using Mellin transform}
\paragraph{}
The Mellin transform of a function \( f(x) \)  is given by:
\[
\mathcal{M}\{f(x)\}(s) = F(s) = \int_0^\infty x^{s-1} f(x) \, dx
\]
where \( s \) is a complex variable.

\begin{theorem}[Mellin inversion formula]
    Followed from the definition of Mellin transform, the Mellin inversion formula is derived to recover the original function \( f(x) \) from its Mellin transform \( F(s) \). The inversion formula is given by:
\[
f(x) = \frac{1}{2\pi i} \int_{c - i\infty}^{c + i\infty} x^{-s} F(s) \, ds
\]
where $c$ is chosen such that the contour path of the complex integration lies within the region of convergence of $F(s)$. 
\end{theorem}
\begin{proof}
See section 2.1.2 of \cite{bertrand1995mellin}. 
\end{proof}

\begin{definition}[Beta function]
    The Beta function, denoted by \( B(x, y) \), is defined by the following integral for two real or complex numbers \( x \) and \( y \) with positive real parts:
\[
B(x, y) = \int_0^1 t^{x-1} (1-t)^{y-1} \, dt
\]
\end{definition}
The next theorem states an important relationship between Beta function and Gamma function. 
\begin{theorem}
    For $x,y\in \mathbb{C}$ where $\Re(x),\Re(y)>0$,  
\[
B(x, y) = \frac{\Gamma(x) \Gamma(y)}{\Gamma(x + y)}
\]
\label{thebeta}
\end{theorem}
\begin{proof}
    See Theorem 0.8 in Chapter \textit{II.0.} of
    \cite{tenenbaum2015introduction}. 
\end{proof}

After introducing some preliminaries, we begin to prove lemma \ref{lemma210}. Mellin transform and its inversion formula derives that if 
$$\mathcal{M}\{f(u)\}(z) = F(z) = \int_{0}^{\infty}{u^{z-1}f(u)du}$$
then, 
\begin{equation}
    \mathcal{M}^{-1}\{F(z)\}(u)=f(u)=\frac{1}{2\pi i} \int_{c- i\infty}^{c + i\infty} F(z) u^{-z} \, dz \label{eq3}
\end{equation}

Hence, we suppose $F(z)$ is the Mellin transform of $f(u)$,
\begin{equation}
    \int_{0}^{\infty}{u^{z-1}f(u)du}=F(z)=\Gamma(z)\Gamma(s-z) \label{eq4}
\end{equation}

Since we have verified the absolute convergence of $\Gamma(z)\Gamma(s-z)$ as integrand in proposition \ref{propsotioncool}. Thus, for any $c$ such that $c\ge \frac{1}{2}$ and $\Re(s)-c\ge \frac{1}{2}$, we have: 
$$f(u)=\frac{1}{2\pi i} \int_{c- i\infty}^{c + i\infty} \Gamma(z) \Gamma(s - z) u^{-z} \, dz$$
and our main goal is to get an explicit expression for $f(u)$.

Theorem \ref{thebeta} gives us: 
$$\frac{\Gamma(z)\Gamma(s-z)}{\Gamma(s)}=B(z,s-z)=B(s-z,z)=\int_0^1 t^{s-z-1} (1-t)^{z-1} \, dt$$

We first substitute $t$ with $\frac{1}{w}$: 
\begin{align*}
    \frac{\Gamma(z)\Gamma(s-z)}{\Gamma(s)}&=\int_{\infty}^1 \left(\frac{1}{w}\right)^{s-z-1} \left(1-\frac{1}{w}\right)^{z-1} \, \left(-\frac{dw}{w^2} \right)\\
    &=\int_{1}^{\infty}{w^{-s}(w-1)^{z-1}dw}
\end{align*}

We then substitute $w$ with $u+1$, which gives: 
$$\Gamma(z)\Gamma(s-z)=\int_0^{\infty}{(u+1)^{-s}\Gamma(s)u^{z-1}}du$$

From \eqref{eq4}: 
$$\Gamma(z)\Gamma(s-z)=\int_0^{\infty}{(u+1)^{-s}\Gamma(s)u^{z-1}}du=\int_{0}^{\infty}{f(u)u^{z-1}du}$$

Therefore, a possible expression of $f(u)$ is $\Gamma(s)(1+u)^{-s}$. A well-known property states that Mellin transform is injective, which means if $\mathcal{M}\{f(x)\}(z)=\mathcal{M}\{g(x)\}(z)$, then $f(x)$ and $g(x)$ must be the same (See Theorem 2.1.2 \cite{bertrand1995mellin}). Therefore, $f(u)$ can only be $\Gamma(s)(1+u)^{-s}$, which finishes the proof. 

\subsubsection{Second proof using Cauchy's residue theorem}
\label{sectiongood}
\paragraph{}
In this method, we avoid the use of Mellin transform. We directly evaluate the complex integral by shifting the integration path and analyzing its isolated singularities.

First, we adjust the linear integration path to a rectangular contour $[c-iT,c+iT], [c+iT,c-k+iT], [c-k+iT,c-k-iT], [c-k-iT,c-iT]$ where $c\ge \frac{1}{2}$ and $\Re(s)-c\ge \frac{1}{2}$. By Cauchy's residue theorem,
\begin{equation}
    \frac{1}{2\pi i}\left[ \int_{c-iT}^{c+iT}+\int_{c+iT}^{c-k+iT}+\int_{c-k+iT}^{c-k-iT}+\int_{c-k-iT}^{c-iT} \right]\Gamma(z)\Gamma(s-z)u^{-z}dz
    \label{eqc}
\end{equation}
equals to the sum of the residues of all isolated singularities within the rectangular contour. 

Next, we need to calculate the sum of the residues. The following theorem gives us an overview to the singularities of Gamma function:
\begin{theorem}
    The function $\Gamma(s)$ initially defined for $\Re(s) > 0$ has an analytic continuation to a meromorphic function on $\mathbb{C}$ whose only singularities are simple poles at the negative integers $s = 0, -1, \ldots$. The residue of $\Gamma(s)$ at $s = -n$ is $\frac{(-1)^n}{n!}$.\label{theorem}
\end{theorem} 
\begin{proof}
    See Theorem 1.3 in Chapter 6 of \cite{stein2003complex}. 
\end{proof}
Hence, within the given rectangular contour where $k$ approaches infinity, the integrand $F(z)=\Gamma(z)\Gamma(s-z)u^{-z}$ can be considered as a meromorphic function with simple poles at all non-positive integers. For each integer $-n$, $n\ge 0$, the residue of $F(z)$ at $z=-n$ can be calculated:
$$\text{Res}\left[\Gamma(z)\Gamma(s-z)u^{-z}, z=-n \right]=\frac{(-1)^n}{n!}\Gamma(s+n)u^{n}$$

From \eqref{eqc} we can get: 
\begin{equation}
    \frac{1}{2\pi i}\left[ \int_{c-iT}^{c+iT}+\int_{c+iT}^{c-k+iT}+\int_{c-k+iT}^{c-k-iT}+\int_{c-k-iT}^{c-iT} \right]\Gamma(z)\Gamma(s-z)u^{-z}dz=\sum_{n=0}^{\infty}{\frac{(-1)^n}{n!}\Gamma(s+n)u^{n}} \label{eqz}
\end{equation} 

To the vertical integration at negative infinity $[c-k+iT,c-k-iT]$, as $T$ approaches infinity, proposition \ref{propositioncoo} states that: 
\begin{equation}
    \lim_{k, T\to \infty}{\frac{1}{2\pi i}}\int_{c-k+iT}^{c-k-iT}\Gamma(z)\Gamma(s-z)u^{-z}dz=0 \label{eqw}
\end{equation}

To the horizontal integration $[c+iT,c-k+iT]$, as $T$ approaches infinity, lemma \ref{con:8} states that: 
\begin{align*}
\left| \int_{c+iT}^{c-k+iT}\Gamma(z)\Gamma(s-z)u^{-z}dz  \right| &\le 
\int_{c+iT}^{c-k+iT}\left| \Gamma(z)\Gamma(s-z)u^{-z} \right|\cdot|dz|\\
&\le \int_{c}^{c-k} C_1 e^{-T} C_2(s) e^{-|\Im(s)-T|} u^{-\sigma} d\sigma \\
&\le \int_{c}^{c-k} C_1 e^{-T} C_2(s) e^{-T+\Im(s)} u^{-\sigma} d\sigma \ \ \ \ \text{as}\  T\to \infty \\
&\to 0 
\end{align*}

Thus, 
\begin{equation}
    \lim_{k, T\to \infty}{\frac{1}{2\pi i}}\int_{c+iT}^{c-k-iT}\Gamma(z)\Gamma(s-z)u^{-z}dz=0 \label{eqx}
\end{equation}

Similarly, consider the horizontal integration $[c-k-iT,c-iT]$, as $T$ approaches infinity, 
\begin{equation}
    \lim_{k, T\to \infty}{\frac{1}{2\pi i}}\int_{c-k-iT}^{c-iT}\Gamma(z)\Gamma(s-z)u^{-z}dz=0  \label{eqy}
\end{equation}

Combine the above results, we derive
\begin{equation}
    \frac{1}{2\pi i} \int_{c- i\infty}^{c + i\infty} \Gamma(z) \Gamma(s - z) u^{-z} \, dz=\sum_{k=0}^{\infty}{\frac{\Gamma(s+k)}{k!}(-u)^k} \label{eqcool}
\end{equation}

Newton's generalized binomial theorem gives that, 
$$(1+u)^{-s}=\sum_{k=0}^{\infty}{\binom{-s}{k} u^k}$$

Thus, 
\begin{align*}
\Gamma(s)(1 + u)^{-s} &= \sum_{k=0}^{\infty}{\binom{-s}{k}\Gamma(s) u^k}\\
&= \sum_{k=0}^{\infty}{\frac{(-s)(-s-1)...(-s-k+1)}{k!} \Gamma(s) u^k}\\
&= \sum_{k=0}^{\infty}{(-1)^k\cdot \frac{s(s+1)...(s+k-1)\Gamma(s)}{k!} u^k}\\
&= \sum_{k=0}^{\infty} (-1)^k \frac{\Gamma(s + k)}{k!} u^k
\end{align*}

Hence, from \eqref{eqcool}, we finish the proof: 
$$ \frac{1}{2\pi i} \int_{c- i\infty}^{c + i\infty} \Gamma(z) \Gamma(s - z) u^{-z} \, dz=\sum_{k=0}^{\infty}{\frac{\Gamma(s+k)}{k!}(-u)^k}=\Gamma(s)(1 + u)^{-s}$$

\subsection{Representation of $\sum_{m=1}^{\infty}\sum_{n=1}^{\infty}{(m+n)^{-s}}$} \label{section23}
\paragraph{}
In this section, we deduce an integral representation of the double sum version of Riemann zeta function $\sum_{m=1}^{\infty}\sum_{n=1}^{\infty}{(m+n)^{-s}}$, which we will further evaluate in the next section. 
\begin{lemma}
For $a, b \in [1, \infty)$, $c\ge \frac{1}{2}$, and $\Re(s)-c\ge \frac{1}{2}$, we have
\begin{equation}
\frac{\Gamma(s)}{(a + b)^s} = \frac{1}{2\pi i} \int_{c - i\infty}^{c + i\infty} \frac{\Gamma(z) \Gamma(s - z)}{a^z b^{s - z}} \, dz   \label{eqcomplex}
\end{equation}
\label{lemmanice}
\end{lemma} 

\begin{proof}
If $a=b$, the result follows naturally by choosing $u=1$ in \eqref{eq2}. Otherwise, without the loss of generality, we assume $a<b$. Thus, we substitute $u$ in \eqref{eq2} by $\frac{a}{b}$ since $\frac{a}{b}\in (0,1)$: 
$$\Gamma(s)\left(1 + \frac{a}{b}\right)^{-s} =\frac{1}{2\pi i} \int_{c- i\infty}^{c + i\infty} \Gamma(z) \Gamma(s - z) \left(\frac{a}{b}\right)^{-z}  dz $$

Hence, 
\begin{align*}
\frac{\Gamma(s)}{(a+b)^s}&= b^{-s}\Gamma(s)\left(1 + \frac{a}{b}\right)^{-s}\\
&= \frac{1}{2\pi i} \int_{c- i\infty}^{c + i\infty} \Gamma(z) \Gamma(s - z) \left(\frac{a}{b}\right)^{-z} b^{-s} \, dz\\
&=\frac{1}{2\pi i} \int_{c - i\infty}^{c + i\infty} \frac{\Gamma(z) \Gamma(s - z)}{a^z b^{s - z}} \, dz 
\end{align*}
\end{proof}

\begin{lemma}
Given complex number $s$ satisfying $\Re(s)>2$, choose $c$ such that $c>1$, and $\Re(s)-c>1$, we have
\[
\Gamma(s) \sum_{m=1}^{\infty} \sum_{n=1}^{\infty} (m + n)^{-s} = \frac{1}{2\pi i} \int_{c - i\infty}^{c + i\infty} \zeta(z) \zeta(s - z) \Gamma(z) \Gamma(s - z) \, dz
\]
where \(\zeta\) denotes the Riemann zeta function.
\end{lemma}

\begin{proof}
The result from lemma \ref{lemmanice} gives that, 
$$\frac{\Gamma(s)}{(m + n)^s} = \frac{1}{2\pi i} \int_{c - i\infty}^{c + i\infty} \frac{\Gamma(z) \Gamma(s - z)}{m^{z} n^{s - z}} \, dz$$

Hence, 
\begin{align*}
    \sum_{m=1}^{\infty} \sum_{n=1}^{\infty}\frac{\Gamma(s)}{(m + n)^s}&=\sum_{m=1}^{\infty} \sum_{n=1}^{\infty}\frac{1}{2\pi i} \int_{c - i\infty}^{c + i\infty} \frac{\Gamma(z) \Gamma(s - z)}{m^{z} n^{s - z}} \, dz\\
    &=\frac{1}{2\pi i} \int_{c - i\infty}^{c + i\infty} \Gamma(z) \Gamma(s - z)\left( \sum_{m=1}^{\infty} \sum_{n=1}^{\infty}\frac{1}{m^{z} n^{s - z}}\right) \, dz\\
    &=\frac{1}{2\pi i} \int_{c - i\infty}^{c + i\infty} \Gamma(z) \Gamma(s - z)\left( \sum_{m=1}^{\infty} \frac{1}{m^{z}}\right)\left( \sum_{n=1}^{\infty} \frac{1}{n^{s - z}}\right) \, dz\\
    &=\frac{1}{2\pi i} \int_{c - i\infty}^{c + i\infty} \zeta(z) \zeta(s - z) \Gamma(z) \Gamma(s - z) \, dz
\end{align*}
which gives the desired result. Lastly, we need to verify the absolute convergence of the integrand, and validate the interchanging of summation and integration within the procedure above. Riemann zeta function $\zeta(s)$ has the abscissa of absolute convergence $\sigma_{a}=1$, which means that the series converges absolutely when the real part of its complex variable $s$ is larger than 1. Thus, when $\Re(z)=c>1$ and $\Re(s)-c>1$, $\zeta(z)\zeta(s-z)$ converges absolutely. Gamma function $\Gamma(s)$ converges absolutely in the half-plane $\Re(s)>0$. Hence, when $\Re(z)=c>1$ and $\Re(s)-c>1$, the integrand $F(z)=\zeta(z)\zeta(s-z)\Gamma(z)\Gamma(s-z)$ converges absolutely. Due to Lebesgue's dominated convergence theorem, it is valid to interchange summation and integration in the above derivation. 
\end{proof}

\subsection{Evaluation of the complex integral} \label{section24}
\paragraph{}
Similar to what we have done previously in section \ref{sectiongood}, we first need to adjust the linear integration path to a rectangular contour $[c-iT,c+iT], [c+iT,c-k+iT], [c-k+iT,c-k-iT], [c-k-iT,c-iT]$ where $c>1$, $\Re(s)-c>1$, and $T, k$ approaches infinity. By Cauchy's residue theorem,
\begin{equation}
    \frac{1}{2\pi i}\left[ \int_{c-iT}^{c+iT}+\int_{c+iT}^{c-k+iT}+\int_{c-k+iT}^{c-k-iT}+\int_{c-k-iT}^{c-iT} \right]\zeta(z) \zeta(s - z) \Gamma(z) \Gamma(s - z) \, dz
    \label{eql}
\end{equation}
equals to the sum of the residues of all simple poles within the rectangular contour. 
\subsubsection{Sum of appropriate residues within the rectangular contour}
\paragraph{}
Next, we need to calculate the sum of the residues. Since we have been acquainted with the analytic properties of Gamma function in the previous proofs, we now need some basic properties of Riemann zeta function. 
\begin{theorem}[Analytic continuation of Riemann zeta function]
    The Riemann zeta function, defined for \( \sigma > 1 \) by the series$$
\zeta(s) = \sum_{n=1}^{\infty} \frac{1}{n^s},$$
has an analytic continuation to a function defined on the half-plane \( \sigma > 0 \) and is analytic in this half-plane with the exception of a simple pole at \( s = 1 \) with residue 1, given by 
\begin{equation}
    \zeta(s) = \frac{s}{s-1} - s \int_1^{\infty} \{x\} x^{-s-1} dx \quad (\sigma > 0)  \label{eqnice}
\end{equation}
\end{theorem}
\begin{proof}
    See Theorem 4.11 in Chapter 4 of \cite{hildebrand2005introduction}.
\end{proof}
\begin{theorem}
    For each $n=1,2,3,...$, $$\zeta(-2n)=0$$
    These are the so-called trivial zeros of Riemann zeta function. 
\end{theorem}
\begin{proof}
    This follows naturally from the functional equation of Riemann zeta function (See Theorem 12.7 in Chapter 12 of \cite{apostol1976introduction}):
    \begin{equation}
        \zeta(1-s)=2(2\pi)^{-s}\Gamma(s)\cos\left(\frac{\pi s}{2} \right)\zeta(s)  \label{eqim}
    \end{equation}
    because $\cos\left(\frac{\pi}{2}(2n+1)\right)=0$ for each positive integer $n$. 
\end{proof}
In order to calculate the sum of residues, we need to analyze the set of simple poles for the integrand $F(z)=\zeta(z)\zeta(s-z)\Gamma(z)\Gamma(s-z)$ within the rectangular contour. Since $c$ is larger than 1, $\zeta(z)$ yields a simple pole at $z=1$ and a set of zeros at all even negative integers. Since $\Gamma(z)$ yields simple poles at all non-positive integers, the only simple poles for $F(z)$ are at 1, 0, and $\{-2m-1\}_{m=0}^{\infty}$. Hence, the residues can be calculated: 
\begin{align*}
\text{Res}\left[ \zeta(z)\zeta(s-z)\Gamma(z)\Gamma(s-z), z=1   \right]&=\lim_{z\to 1}{(z-1)\zeta(z)\zeta(s-z)\Gamma(z)\Gamma(s-z)} \\
&=\zeta(s-1)\Gamma(1)\Gamma(s-1)\\
&=\zeta(s-1)\Gamma(s-1)
\end{align*}
\begin{align*}
\text{Res}\left[ \zeta(z)\zeta(s-z)\Gamma(z)\Gamma(s-z), z=0   \right]&=\lim_{z\to 0}{z\zeta(z)\zeta(s-z)\Gamma(z)\Gamma(s-z)}    \\
&=\lim_{z\to 0}{\zeta(z)\zeta(s-z)\Gamma(z+1)\Gamma(s-z)}\\
&=-\frac{1}{2}\zeta(s)\Gamma(s)
\end{align*}
since $\zeta(0)=-\frac{1}{2}$ can be calculated using \eqref{eqim}: 
\begin{align*}
-1=\lim_{z\to 1}(1-z)\zeta(z)&=\lim_{z\to 1}(1-z) 2^z \pi^{z-1} \sin\left(\frac{\pi z}{2}\right) \Gamma(1-z) \zeta(1-z)    \\
&=\lim_{z\to 1}\Gamma(2-z)2^z\pi^{z-1} \sin\left(\frac{\pi z}{2}\right)\zeta(1-z)\\
&=2\zeta(0)
\end{align*}
\begin{align*}
\text{Res}\left[\zeta(z)\zeta(s-z)\Gamma(z)\Gamma(s-z), z=-2m-1  \right]&=\lim_{z\to -2m-1} (z+2m+1)\zeta(z)\zeta(s-z)\Gamma(z)\Gamma(s-z)\\
&=\lim_{z\to -2m-1} \frac{(z+2m+1)(z+2m)...z\Gamma(z)}{(z+2m)(z+2m-1)...z} \zeta(z)\zeta(s-z)\Gamma(s-z)\\
&=\lim_{z\to -2m-1} \frac{\Gamma(z+2m+2)}{(z+2m)(z+2m-1)...z} \zeta(z)\zeta(s-z)\Gamma(s-z)\\
&=-\frac{1}{(2m+1)!}\zeta(-2m-1)\zeta(s+2m+1)\Gamma(s+2m+1)
\end{align*}

Thus, combine the above residues and \eqref{eql}, we get: 
\begin{equation}
\begin{aligned}
&\frac{1}{2\pi i}\left[ \int_{c-iT}^{c+iT}+\int_{c+iT}^{c-k+iT}+\int_{c-k+iT}^{c-k-iT}+\int_{c-k-iT}^{c-iT} \right] \zeta(z) \zeta(s - z) \Gamma(z) \Gamma(s - z) dz\\
=&\zeta(s-1)\Gamma(s-1)-\frac{1}{2}\zeta(s)\Gamma(s)-\sum_{m=0}^{\infty}\frac{1}{(2m+1)!} \zeta(-2m-1)\zeta(s+2m+1)\Gamma(s+2m+1)
\end{aligned}
\label{eq14}
\end{equation}
\subsubsection{Evaluation of the horizontal integral} \label{section242}
\paragraph{}
In this section, we evaluate the horizontal integral on $[c+iT,c-k+iT]$ and $[c-k-iT,c-iT]$ by using the upper bound of Gamma and Riemann zeta function in each assigned region. We first prove a few lemmas about the bound for zeta function under two different conditions. 
\begin{lemma}[Bound for zeta function on the left half-plane $\sigma<0$]
In any vertical strip on the left half-plane $\sigma<0$, we have: 
$$|\zeta(\sigma + it)| = \mathcal{O}\left(|t|^{\frac{1}{2}-\sigma}\right)$$
\end{lemma}
\begin{proof}
The functional equation of Riemann zeta function \eqref{eqim} derives: 
$$\zeta(s) = 2^s \pi^{s-1} \sin \left( \frac{\pi s}{2} \right) \Gamma(1-s) \zeta(1-s)$$

Since $\sigma<0$, we have
$$|\zeta(1-\sigma-it)|=\left| \sum_{n=1}^{\infty}\frac{1}{n^{1-\sigma-it}} \right|\le \sum_{n=1}^{\infty}\frac{1}{n^{1-\sigma}}=\zeta(1-\sigma) $$
converges to a fixed value. Hence, combine the above results and Lemma:
$$|\zeta(\sigma+it)|=\mathcal{O}\left(e^{\frac{\pi}{2}|t|}|t|^{1-\frac{1}{2}-\sigma}e^{-\frac{\pi}{2}|t|}\right)=\mathcal{O}\left(|t|^{\frac{1}{2}-\sigma}\right)$$
\end{proof}
\begin{lemma}[Bound for zeta function in the vertical strip $0<\sigma<1$]
Given $0<\delta<1$, for $\delta \le \sigma<1$ and $|t|\ge 2$, we have: 
$$|\zeta(\sigma + it)|= \mathcal{O}(|t|^{1-\delta}) $$
\end{lemma}
\begin{proof}
Lemma 5.4 in Chapter 5 of \cite{hildebrand2005introduction} states: 
$$\zeta(s) = \sum_{n=1}^{N} \frac{1}{n^s} - \frac{N^{1-s}}{1-s} - s \int_{N}^{\infty} \frac{\{x\}}{x^{s+1}} dx $$ for some integer $N$.
Hence, for $\sigma\le \delta$, we have: 
\begin{align*}
|\zeta(\sigma+it)|&\le \left|\sum_{n=1}^{N} \frac{1}{n^s}  \right| + \left|\frac{N^{1-s}}{1-s}\right| + \left| s \int_{N}^{\infty} \frac{\{x\}}{x^{s+1}} dx \right| \\ 
 & =\sum_{n=1}^{N}{\frac{1}{n^{\delta}}}+\frac{N^{1-\delta}}{1-\delta}+|\sigma+it|\int_{N}^{\infty}{x^{-\delta-1}dx}\\
 &=\mathcal{O}(N^{1-\delta})+\frac{|\sigma+it|}{\delta N^{\delta}}
\end{align*}

We choose $N=\lfloor |t| \rfloor$, then it gives the desired result: 
$$|\zeta(\sigma + it)|= \mathcal{O}(|t|^{1-\delta})$$
\end{proof}

After introducing all necessary results, we first evaluate the horizontal integral on $[c+iT,c-k+iT]$. We spilt the integral into two segments in terms of the range of the variable in Riemann zeta function. 
\begin{equation}
\begin{aligned}
    & \ \ \   \left| \frac{1}{2\pi i}\int_{c+iT}^{c-k+iT} \zeta(z) \zeta(s - z) \Gamma(z) \Gamma(s - z) \, dz \right|  \\
    &\le \frac{1}{2\pi} \int_{c-k}^{c} \left| \zeta(\sigma+iT) \zeta(s-\sigma-iT) \Gamma(\sigma+iT) \Gamma(s-\sigma-iT) \right| d\sigma \\
    &=\frac{1}{2\pi} \left[ \int_{c-k}^{0}+\int_{0}^{c} \right] \left| \zeta(\sigma+iT) \zeta(s-\sigma-iT) \Gamma(\sigma+iT) \Gamma(s-\sigma-iT) \right| d\sigma  \\
    &\ll \frac{1}{2\pi}\left[\int_{c-k}^{0} T^{\frac{1}{2}-\sigma}\zeta(s) C_1 e^{-T} d\sigma + \int_{0}^{c} T\zeta(s-1) C_2(s) e^{-T} d\sigma \right]   \\ 
    &=\frac{1}{2\pi}\left[\frac{T^{\frac{1}{2}-c}(T^k-T^c)}{\log(T)} C_1 \zeta(s)e^{-T}+cT\zeta(s-1)C_{2}(s)e^{-T}\right] 
   \label{eqdi}
\end{aligned}
\end{equation}

As $T$ approaches infinity, equation \eqref{eqdi} tends to 0. Hence,
$$\lim_{T\to \infty}{\frac{1}{2\pi i}\int_{c+iT}^{c-k+iT} \zeta(z) \zeta(s - z) \Gamma(z) \Gamma(s - z) \, dz}=0$$

Similarly, on the horizontal segment $[c-k-iT,c-iT]$, we can derive: 
$$\lim_{T\to \infty}{\frac{1}{2\pi i}\int_{c-k-iT}^{c-iT} \zeta(z) \zeta(s - z) \Gamma(z) \Gamma(s - z) \, dz}=0$$

\subsubsection{Evaluation of the vertical integral}
\paragraph{}
Finally, we need to evaluate the integral on the vertical segment $[c-k+iT,c-k-iT]$. Before handling this integral, we introduce an important limit regarding Riemann zeta function. 
\begin{lemma}
Suppose $\sigma=\Re(s)$, then we have: 
\[
\lim_{\sigma \to \infty} \zeta(s)=1
\label{lemmalimit}
\]
\end{lemma}
\begin{proof}
Due to the definition of Riemann zeta function when the real part of $s$ is larger than 1: 
$$\lim_{\sigma\to \infty} \zeta(s)=\lim_{\sigma\to \infty} \sum_{n=1}^{\infty} \frac{1}{n^s}$$

Tannery's theorem \cite{tannery_wikipedia} allows us to interchange the limit and summation. Hence, 
$$\lim_{\sigma\to \infty} \sum_{n=1}^{\infty} \frac{1}{n^s}=\sum_{n=1}^{\infty} \lim_{\sigma\to \infty} \frac{1}{n^s}$$

Since
$$\lim_{\sigma\to \infty}\left|\frac{1}{n^s}\right|=\lim_{\sigma\to \infty}\frac{1}{n^\sigma}=0 \ \ \ \ (n>1)$$
then, 
\[
\lim_{\sigma \to \infty} \zeta(s)=1
\]
\end{proof}

Then, we solve the integral by first substituting the variable $z$ by $z-k$:
\[
\frac{1}{2\pi i}\int_{c-k+iT}^{c-k-iT}\zeta(z) \zeta(s - z) \Gamma(z) \Gamma(s - z) \, dz=\frac{1}{2\pi i}\int_{c+iT}^{c-iT}\zeta(z-k) \zeta(s - z+k) \Gamma(z-k) \Gamma(s - z+k)  dz
\]

Using the functional equation of Riemann zeta function: 
\[
=\frac{1}{2\pi i}\int_{c+iT}^{c-iT} 2^{z-k}(\pi)^{z-k-1}\sin\left(\frac{\pi}{2}(z-k)\right)\zeta(1-z+k)\zeta(s-z+k)\Gamma(1-z+k)\Gamma(z-k)\Gamma(s-z+k)dz
\]

We use the Euler's reflection formula $\Gamma(z)\Gamma(1-z)=\frac{\pi}{\sin(\pi z)}$ to substitute sine function by Gamma function: 
\[=\frac{1}{2\pi i}\int_{c+iT}^{c-iT} (2\pi)^{z-k}\zeta(1-z+k)\zeta(s-z+k) \ \frac{\Gamma(1-z+k)\Gamma(z-k)\Gamma(s-z+k)}{\Gamma(\frac{z}{2}-\frac{k}{2}) \Gamma(1-\frac{z}{2}+\frac{k}{2})} dz
\]

We then use Legendre's duplication formula $\Gamma(2z)=2^{2z-1}\pi^{-\frac{1}{2}}\Gamma(z)\Gamma(z+\frac{1}{2})$ to substitute $\Gamma(1-z+k)$ and $\Gamma(z-k)$ in the above equation: 
\begin{align*}
&=\frac{1}{2\pi i}\int_{c+iT}^{c-iT} (2\pi)^{z-k}\zeta(1-z+k)\zeta(s-z+k)\frac{\Gamma(\frac{1}{2}-\frac{z}{2}+\frac{k}{2})\Gamma(1-\frac{z}{2}+\frac{k}{2}) \Gamma(\frac{z}{2}-\frac{k}{2})\Gamma(\frac{1}{2}+\frac{z}{2}-\frac{k}{2})\Gamma(s-z+k)}{2\pi \Gamma(\frac{z}{2}-\frac{k}{2})\Gamma(1-\frac{z}{2}+\frac{k}{2})}dz\\
&=\frac{1}{2\pi i}\int_{c+iT}^{c-iT} (2\pi)^{z-k-1} \zeta(1-z+k)\zeta(s-z+k) \Gamma\left(\frac{1}{2}-\frac{z}{2}+\frac{k}{2}\right) \Gamma\left(\frac{1}{2}+\frac{z}{2}-\frac{k}{2}\right) \Gamma(s-z+k) dz
\end{align*}

Using lemma \ref{lemmalimit} $\lim_{\Re(s)\to \infty} \zeta(s)=1$, we can simplify the expression as $k$ approaches infinity: 
\begin{equation}
=\lim_{k \to \infty}\frac{1}{2\pi i}\int_{c+iT}^{c-iT} (2\pi)^{z-k-1} \Gamma\left(\frac{1}{2}-\frac{z}{2}+\frac{k}{2}\right) \Gamma\left(\frac{1}{2}+\frac{z}{2}-\frac{k}{2}\right) \Gamma(s-z+k) dz  
\label{eqsolve}
\end{equation}

Hence, our next step is to solve this expression. In order to handle it, we introduce Hurwitz zeta function, which is defined as: 
$$\zeta(s,a)=\sum_{n=0}^{\infty}{\frac{1}{(n+a)^s}}$$
for $\Re(s)>1$ and $a\neq 0,1,2,...$. In \cite{apostol1976introduction}, T.M.Apostol expresses Hurwitz zeta function using periodic zeta function $\sum_{n=1}^{\infty}{e^{2\pi ina}/ n^s}$: 
$$\zeta(s,a)=\frac{\Gamma(1-s)}{(2\pi)^{1-s}}\left[ e^{-\pi i(1-s)/2} \sum_{n=1}^{\infty}{\frac{e^{2\pi ina}}{n^{1-s}}}+e^{\pi i(1-s)/2} \sum_{n=1}^{\infty}{\frac{e^{-2\pi ina}}{n^{1-s}}} \right]$$

Alternatively, we can represent Hurwitz zeta function by a complex integral using lemma \ref{lemmanice}: 
$$\frac{\Gamma(s)}{(n+(a-1))^s}=\frac{1}{2\pi i}\int_{c-i\infty}^{c+i\infty} \frac{\Gamma(z)\Gamma(s-z)}{n^z (a-1)^{s-z}}dz$$
for some $n\ge 1$ and $a\ge 2$. Summing $n$ from 1 to infinity: 
$$\Gamma(s)\zeta(s,a)=\sum_{n=0}^{\infty} \frac{\Gamma(s)}{(n+a)^s}=\sum_{n=1}^{\infty} \frac{\Gamma(s)}{(n+(a-1))^s}=\frac{1}{2\pi i}\sum_{n=1}^{\infty}\int_{c-i\infty}^{c+i\infty} \frac{\Gamma(z)\Gamma(s-z)}{n^z (a-1)^{s-z}}dz$$

Lebesgue's dominated convergence theorem allows us to interchange the summation and integration:
\begin{align*}
\Gamma(s)\zeta(s,a)&=\frac{1}{2\pi i} \int_{c-i\infty}^{c+i\infty} \sum_{n=1}^{\infty} \frac{\Gamma(z)\Gamma(s-z)}{n^z (a-1)^{s-z}}dz\\
&=\frac{1}{2\pi i} \int_{c-i\infty}^{c+i\infty} \zeta(z)\Gamma(z)\Gamma(s-z)(a-1)^{-s+z}dz 
\end{align*}

On one hand, this complex integral equals to $\Gamma(s)\zeta(s,a)$. On the other hand, we can shift the integration path to a rectangular contour to evaluate it. Same to what we have done previously, we choose the rectangle $[c-iT,c+iT], [c+iT,c-k+iT], [c-k+iT,c-k-iT], [c-k-iT,c-iT]$ where $c>1$, $\Re(s)-c>1$, and $T, k$ approaches infinity. By Cauchy's residue theorem,
\begin{equation}
    \frac{1}{2\pi i}\left[ \int_{c-iT}^{c+iT}+\int_{c+iT}^{c-k+iT}+\int_{c-k+iT}^{c-k-iT}+\int_{c-k-iT}^{c-iT} \right]\zeta(z) \Gamma(z) \Gamma(s - z)(a-1)^{-s+z} \, dz
    \label{eqlmi}
\end{equation}
equals to the sum of the residues of all isolated singularities within the rectangular contour. 
\begin{align*}
\text{Res}[\zeta(z)\Gamma(z)\Gamma(s-z)(a-1)^{-s+z},z=1]&=\lim_{z\to 1} (z-1)\zeta(z)\Gamma(z)\Gamma(s-z)(a-1)^{-s+z}\\
&=\Gamma(1)\Gamma(s-1)(a-1)^{-s+1}\\
&=\Gamma(s-1)(a-1)^{-s+1}
\end{align*}
\begin{align*}
\text{Res}[\zeta(z)\Gamma(z)\Gamma(s-z)(a-1)^{-s+z},z=0]&=\lim_{z\to 0} z\zeta(z)\Gamma(z)\Gamma(s-z)(a-1)^{-s+z}    \\
&=\lim_{z\to 0} \zeta(z)\Gamma(z+1)\Gamma(s-z)(a-1)^{-s+z}\\
&=-\frac{1}{2}\Gamma(s)(a-1)^{-s}
\end{align*}
\begin{align*}
\text{Res}[\zeta(z)\Gamma(z)\Gamma(s-z)(a-1)^{-s+z},z=-2m-1]&=\lim_{z\to -2m-1} (z+2m+1)\zeta(z)\Gamma(z)\Gamma(s-z)(a-1)^{-s+z} \\
&=\lim_{z\to -2m-1} \frac{(z+2m+1)...z\Gamma(z)}{(z+2m)...z}\zeta(z)\Gamma(s-z)(a-1)^{-s+z}\\ 
&=\lim_{z\to -2m-1} \frac{\Gamma(z+2m+2)}{(z+2m)...z}\zeta(z)\Gamma(s-z)(a-1)^{-s+z}\\
&=-\frac{1}{(2m+1)!}\zeta(-2m-1)\Gamma(s+2m+1)(a-1)^{-s-2m-1}
\end{align*}

Hence, from \eqref{eqlmi} we get:
\begin{equation}
\begin{aligned}
&\frac{1}{2\pi i}\left[ \int_{c-iT}^{c+iT}+\int_{c+iT}^{c-k+iT}+\int_{c-k+iT}^{c-k-iT}+\int_{c-k-iT}^{c-iT} \right] \zeta(z) \Gamma(z) \Gamma(s - z)(a-1)^{-s+z}dz\\ 
=&\Gamma(s-1)(a-1)^{-s+1}-\frac{1}{2}\Gamma(s)(a-1)^{-s}-\sum_{m=0}^{\infty}\frac{1}{(2m+1)!} \zeta(-2m-1)\Gamma(s+2m+1)(a-1)^{-s-2m-1} 
\end{aligned}
\end{equation}

Using the same bound for Gamma and Riemann zeta function in section \ref{section242}, we can show that as $T$ approaches infinity, the integrals on the horizontal segments tend to zero. Thus, the integral along the vertical segment at negative infinity can be represented as: 
\begin{equation}
\begin{split}
\lim_{k,T\to \infty}\int_{c-k+iT}^{c-k-iT} \zeta(z)\Gamma(z)\Gamma(s-z)(a-1)^{-s+z} =\Gamma(s-1)(a-1)^{-s+1}-\frac{1}{2}\Gamma(s)(a-1)^{-s}\\
-\sum_{m=0}^{\infty}\frac{1}{(2m+1)!} \zeta(-2m-1)\Gamma(s+2m+1)(a-1)^{-s-2m-1}-\Gamma(s)\zeta(s,a)  
\end{split}
\label{eqimpor}
\end{equation} 

Next, we rewrite the integral at negative infinity by substituting $z$ in the integrand by $z-k$:
$$\int_{c-k+iT}^{c-k-iT}=\int_{c+iT}^{c-iT} \zeta(z-k)\Gamma(z-k)\Gamma(s-z+k)(a-1)^{-s+z-k}dz$$

Using the functional equation of Riemann zeta function: 
$$=\int_{c+iT}^{c-iT} 2^{z-k} \pi^{z-k-1}\sin\left(\frac{\pi}{2}(z-k) \right)\zeta(1-z+k)\Gamma(1-z+k)\Gamma(z-k)\Gamma(s-z+k)(a-1)^{-s+z-k} dz$$

As $k$ approaches infinity, we use the limit $\lim_{\Re(s) \to \infty} \zeta(s)=1$ to simplify the equation: 
$$=\int_{c+iT}^{c-iT} 2^{z-k} \pi^{z-k-1}\sin\left(\frac{\pi}{2}(z-k) \right)\Gamma(1-z+k)\Gamma(z-k)\Gamma(s-z+k)(a-1)^{-s+z-k} dz$$

Use Euler's reflection formula and Legendre's duplication formula, we further simplify the result: 
$$=\int_{c+iT}^{c-iT} (2\pi)^{z-k-1}\Gamma\left(\frac{1}{2}-\frac{z}{2}+\frac{k}{2}\right) \Gamma\left(\frac{1}{2}+\frac{z}{2}-\frac{k}{2}\right) \Gamma(s-z+k)(a-1)^{-s+z-k} dz $$

Combine equation \eqref{eqimpor}, we derive: 
\begin{equation}
\begin{aligned}
&\lim_{k,T\to \infty}\int_{c+iT}^{c-iT} (2\pi)^{z-k-1}\Gamma\left(\frac{1}{2}-\frac{z}{2}+\frac{k}{2}\right) \Gamma\left(\frac{1}{2}+\frac{z}{2}-\frac{k}{2}\right) \Gamma(s-z+k)(a-1)^{-s+z-k} dz \\
=&-\frac{1}{2}\Gamma(s)(a-1)^{-s}-\Gamma(s)\zeta(s,a)-\sum_{m=0}^{\infty}\frac{1}{(2m+1)!} \zeta(-2m-1)\Gamma(s+2m+1)(a-1)^{-s-2m-1}
\end{aligned}
\end{equation}
which is what we need to solve in equation \eqref{eqsolve}. Choose $a=2$: 
\begin{align*}
&\lim_{k,T \to \infty}\frac{1}{2\pi i}\int_{c+iT}^{c-iT} (2\pi)^{z-k-1} \Gamma\left(\frac{1}{2}-\frac{z}{2}+\frac{k}{2}\right) \Gamma\left(\frac{1}{2}+\frac{z}{2}-\frac{k}{2}\right) \Gamma(s-z+k) dz \\
=&\Gamma(s-1)-\frac{1}{2}\Gamma(s)-\sum_{m=0}^{\infty}\frac{1}{(2m+1)!}\zeta(-2m-1)\Gamma(s+2m+1)-\Gamma(s)\zeta(s,2)
\end{align*}
which gives us the evaluation of complex integral of $\zeta(z)\zeta(s-z)\Gamma(z)\Gamma(s-z)$ at negative infinity.

\subsubsection{Finalize the evaluation}
\paragraph{}
Equation \eqref{eq14} gives us: 
\begin{equation*}
\begin{split}
\frac{1}{2\pi i}\left[ \int_{c-iT}^{c+iT}+\int_{c+iT}^{c-k+iT}+\int_{c-k+iT}^{c-k-iT}+\int_{c-k-iT}^{c-iT} \right] \zeta(z) \zeta(s - z) \Gamma(z) \Gamma(s - z) \, dz\\
=\zeta(s-1)\Gamma(s-1)-\frac{1}{2}\zeta(s)\Gamma(s)-\sum_{m=0}^{\infty}\frac{1}{(2m+1)!} \zeta(-2m-1)\zeta(s+2m+1)\Gamma(s+2m+1)
\end{split}
\end{equation*}

Hence, as $T$ approaches infinity: 
\begin{equation}
\begin{aligned}
\frac{1}{2\pi i}\int_{c-i\infty}^{c+i\infty}&=\lim_{k\to \infty}\left[\frac{1}{2\pi i}\int_{c-k-i\infty}^{c-k+i\infty}\right]+\zeta(s-1)\Gamma(s-1)-\frac{1}{2}\zeta(s)\Gamma(s)  \\
&\ \ \ -\sum_{m=0}^{\infty}\frac{1}{(2m+1)!} \zeta(-2m-1) \zeta(s+2m+1) 
\Gamma(s+2m+1)  \\
&=-\Gamma(s-1)+\frac{1}{2}\Gamma(s)+\sum_{m=0}^{\infty}\frac{1}{(2m+1)!} \zeta(-2m-1)\Gamma(s+2m+1)+\Gamma(s)\zeta(s,2)+\zeta(s-1)\Gamma(s-1) \\
&\ \ \ -\frac{1}{2}\zeta(s)\Gamma(s)-\sum_{m=0}^{\infty}\frac{1}{(2m+1)!} \zeta(-2m-1) \zeta(s+2m+1) \Gamma(s+2m+1) \\
&=\frac{1}{2}\Gamma(s)\zeta(s,2)+\Gamma(s-1)\zeta(s-1,2)-\sum_{m=0}^{\infty}\frac{1}{(2m+1)!} \zeta(-2m-1) \zeta(s+2m+1,2) \Gamma(s+2m+1) \label{eq29}
\end{aligned}
\end{equation}

\section{Application and further discussion} \label{section3}
In the last section, we discuss how the result from section \ref{section2} can be connected with some other integrals. We first derive an integral representation of power function using the definition of Gamma function: 
$$x^{-s}=\frac{1}{\Gamma(s)}\int_{0}^{\infty} e^{-\lambda} \lambda^{s-1} x^{-s} d\lambda$$

We substitute $\lambda$ with $xt$:
\begin{align*}
x^{-s}&=\frac{1}{\Gamma(s)} \int_{0}^{\infty} e^{-xt}(xt)^{s-1}x^{-s}xdt\\
&=\frac{1}{\Gamma(s)} \int_{0}^{\infty} e^{-xt} t^{s-1} dt
\end{align*}

Hence, the double sum in section \ref{section23} can be rewritten in another integral form: 
$$\sum_{m=1}^{\infty}\sum_{n=1}^{\infty}(m+n)^{-s}=\sum_{m=1}^{\infty}\sum_{n=1}^{\infty}\frac{1}{\Gamma(s)} \int_{0}^{\infty} e^{-(m+n)t} t^{s-1} dt$$

Lebesgue's dominated convergence theorem allows to interchange the summation and integration:
\begin{align*}
&=\frac{1}{\Gamma(s)}\int_{0}^{\infty} \sum_{m=1}^{\infty}\sum_{n=1}^{\infty} e^{-(m+n)t}t^{s-1}dt\\ 
&=\frac{1}{\Gamma(s)} \int_{0}^{\infty} \frac{t^{s-1}}{(e^t-1)^2}dt 
\end{align*}

Therefore,
$$\int_{0}^{\infty} \frac{t^{s-1}}{(e^t-1)^2}dt=\Gamma(s)\sum_{m=1}^{\infty}\sum_{n=1}^{\infty}(m+n)^{-s} $$

Equation \eqref{eq29} gives us the representation of this integral when $s$ is complex number whose real part is larger than 1:
$$\int_{0}^{\infty} \frac{t^{s-1}}{(e^t-1)^2}dt=\frac{1}{2}\Gamma(s)\zeta(s,2)+\Gamma(s-1)\zeta(s-1,2)-\sum_{m=0}^{\infty}\frac{1}{(2m+1)!} \zeta(-2m-1) \zeta(s+2m+1,2) \Gamma(s+2m+1)$$

Alternatively, equation \eqref{eq29} can be expressed in another integral form. Similar to the product formula of Gamma and Riemann zeta function \eqref{cond}, a generalized formula of Gamma and Hurwitz zeta function is famously known as: 
$$\Gamma(s)\zeta(s,a)=\int_{0}^{\infty} \frac{x^{s-1}e^{-ax}}{1-e^{-x}} dx$$

Hence,
$$\Gamma(s)\zeta(s,2)=\int_{0}^{\infty} \frac{x^{s-1}e^{-2x}}{1-e^{-x}} dx=\int_{0}^{\infty} \frac{x^{s-1}}{e^x(e^x-1)}dx$$
$$\Gamma(s-1)\zeta(s-1,2)=\int_{0}^{\infty} \frac{x^{s-2}}{e^x(e^x-1)}dx$$

Substitute $\Gamma(s)\zeta(s,2)$ and $\Gamma(s-1)\zeta(s-1,2)$ by their integrals in equation \eqref{eq29}. Also, we use the representation of Riemann zeta function at negative order in terms of Bernoulli number \cite{edwards1974riemann}: $\zeta(-n)=(-1)^n \frac{B_{n+1}}{n+1}$ where $B_{n+1}$ is the $n+1$-th Bernoulli number. We can then get: 
\begin{align*}
\frac{1}{2\pi i}\int_{c-i\infty}^{c+i\infty}&=\int_{0}^{\infty} \frac{x^{s-1}}{2e^x(e^x-1)}dx+\int_{0}^{\infty}\frac{x^{s-2}}{e^x(e^x-1)}dx+ \sum_{m=0}^{\infty} \frac{B_{2m+2}}{(2m+2)!} \int_{0}^{\infty} \frac{x^{s+2m}}{e^x(e^x-1)}dx\\
&=\int_{0}^{\infty} \frac{x^{s-1}}{2e^x(e^x-1)}dx+\int_{0}^{\infty}\frac{x^{s-2}}{e^x(e^x-1)}dx+\int_{0}^{\infty} \left( \sum_{m=0}^{\infty} \frac{B_{2m+2}\cdot x^{2m+2}}{(2m+2)!} \right) \frac{x^{s-2}}{e^x(e^x-1)} dx\\
&=\int_{0}^{\infty} \frac{x^{s-1}+2x^{s-2}}{2e^x(e^x-1)}dx+\int_{0}^{\infty} \left( \sum_{m=1}^{\infty} \frac{B_{2m}\cdot x^{2m}}{(2m)!} \right) \frac{x^{s-2}}{e^x(e^x-1)} dx\\
&=\int_{0}^{\infty} \frac{x^{s-1}+2x^{s-2}+(x\coth(\frac{x}{2})-2)x^{s-2}}{2e^x(e^x-1)} dx\\
&=\int_{0}^{\infty} \frac{x^{s-1}(\coth(\frac{x}{2})+1)}{2e^x(e^x-1)} dx
\end{align*}

Here, Lebesgue's dominated convergence theorem validates the interchange of summation and integration. We also use the Taylor expansion of hyperbolic cotangent function: 
$$\frac{x}{2}\coth\left(\frac{x}{2}\right)=\sum_{n=0}^{\infty} \frac{B_{2n} x^{2n}}{(2n)!}=B_{0}+\sum_{n=1}^{\infty} \frac{B_{2n} x^{2n}}{(2n)!}$$
where $B_{0}$, according to the definition of Bernoulli number, equals to 1.

In summary, the evaluation of the complex integral \eqref{eq123} establishes a broad and profound connection with integrals involving various special functions, such as the exponential and hyperbolic functions. These connections reveal underlying structural similarities in complex analysis and provide valuable insights that could deepen our understanding of such functions in different mathematical contexts. Further research could focus on extending this analysis to more generalized forms, particularly by exploring the double sum version of the Riemann zeta function:

\[
\sum_{m=1}^{\infty} \sum_{n=1}^{\infty} \frac{1}{(m+n)^s}
\]

This form could be extended to a triple sum, and subsequently to sums of arbitrary dimension. Such generalizations could lead to a richer theory encompassing multi-dimensional summations and their connections with known functions in number theory and complex analysis.

Additionally, a worthwhile avenue for future research lies in exploring the explicit expression of the Barnes zeta function, which represents an extended and more complex version of the traditional Riemann zeta function $\zeta(s)$. The Barnes zeta function is defined as:

\[
\zeta_{N}(s,w|a_1,a_2,\dots,a_N) = \sum_{n_1=0}^{\infty} \dots \sum_{n_N=0}^{\infty} \frac{1}{(w + a_1 n_1 + \dots + a_N n_N)^s}
\]
where \( w \) and \( a_i \) have positive real parts, and \( s \) has a real part greater than \( N \). This function introduces an additional layer of complexity by incorporating multiple parameters and summations, which allows for a vast range of potential applications in fields such as analytic number theory, physics, and mathematical physics. The study of the Barnes zeta function could provide insights into multi-dimensional zeta functions, which can be valuable in leading to new discoveries.

\newpage
\section*{Acknowledgment}
\paragraph{}
Firstly, I want to express my sincere gratitude to Professor Ting Wu from Nanjing University for her free guidance and consistent encouragement throughout this research project. She has provided me with necessary books and articles in mathematics, tutored me in academic writing, and gave me valuable suggestions for revision on my paper. Her consistent encouragement has played a pivotal role in promoting me to think deep and finish this research. Secondly, I want to thank Chinese Association of Science and Technology for initiating China Talent Program which provides fantastic opportunities for high school students to learn advanced topics in mathematics and do original research. Lastly and most importantly, I want to thank my parents who give me countless support not only during this research, but also throughout my entire course of learning mathematics. \\
\newline 

The motivation for this research project originates from a challenging exercise (with no provided solutions) found in Chapter \textit{II.0} of G. Tenenbaum's monograph \textit{Introduction to Analytic and Probabilistic Number Theory} \cite{tenenbaum2015introduction}. In this exercise, an integral involves both the Gamma and Riemann zeta functions is introduced, which is similar to the following expression:

\begin{equation}
\frac{1}{2\pi i}\int_{c-i\infty}^{c+i\infty} \zeta(z)\zeta(s-z)\Gamma(z)\Gamma(s-z)dz
\label{eqla}
\end{equation}

This integral expression poses a unique mathematical challenge as it combines the properties of both the Riemann zeta function and the Gamma function within a complex contour integration framework. However, despite the introduction of such an intricate integral, Tenenbaum's work provides no solutions or hints regarding its evaluation, leaving a gap that this research aims to address.

This paper begins by thoroughly investigating preliminary knowledge regarding the analytic properties of Gamma and Riemann zeta functions, which validates the conditions under which this integral is defined and converges. This analysis lays the groundwork for developing the approach to deriving an explicit expression for the integral. In pursuing this goal, I read through related research papers on Riemann zeta function from arXiv, ResearchGate, and ScienceDirect. When facing challenges in understanding certain techniques or concepts, I searched for questions on Mathematics Stack Exchange and joined online discussions. Eventually, I got inspiration from existing strategies in the literature, particularly from the works of Chavan et al. \cite{chavan2022hurwitz} and Chaudhry et al. \cite{chaudhry2001extended}. These studies utilize the Cauchy residue theorem, Cauchy modular transform, and other zeta functions to evaluate integrals with similar structures. Building on these techniques, I then introduce the Hurwitz zeta function, which serves as a crucial element in simplifying and ultimately solving the integral in \eqref{eqla}. This approach leads to the development of a unique method, distinct from prior techniques, which provides new insights and yields a complete solution to the integral. In the final section of this paper, I explore potential applications of this result. This newfound solution can fill a gap in Tenenbaum's monograph as well as provide directions for further research in analytic number theory, particularly in areas that require evaluations of integrals involving special functions like the Gamma and Riemann zeta functions.

\end{document}